\documentclass[12pt]{amsart}
\usepackage[mathcal]{eucal}
\usepackage{amssymb, graphics, epsfig}

\newtheorem{thm}{Theorem}
\newtheorem{prop}[thm]{Proposition}
\newtheorem{lem}[thm]{Lemma}

\newtheorem{conj}[thm]{Conjecture}

\theoremstyle{remark}
\newtheorem{rem}[thm]{Remark}

\theoremstyle{definition}

\newcommand{\col}{\kern -3pt :}
\newcommand{\C}{\mathbb C}

\newcommand{\Id}{\mathrm{Id}}
\newcommand{\End}{\mathrm{End}}
\newcommand{\Hom}{\mathrm{Hom}}

\newcommand{\T}{\mathcal T}
\newcommand{\W}{\mathcal W}
\newcommand{\hsp}{\hspace{3 mm}}

\renewcommand{\leq}{\leqslant}

\renewcommand{\phi}{\varphi}

\def\ignore#1{\relax}

\title[QTS and $6j$-symbol]
{Quantum Teichm\"uller spaces\\ and Kashaev's $6j$-symbols}
\author{Hua Bai}
\address{Department
of Mathematics, University of Georgia, Athens GA, U.S.A.}
 \email{huabai@math.uga.edu}
\urladdr{http://www.math.uga.edu/\~{}%
huabai}

\date{\today}

\begin{document}

\begin{abstract}

The Kashaev invariants of $3$-manifolds are based on $6j$-symbols
from the representation theory of the Weyl algebra, a Hopf algebra
corresponding to the Borel subalgebra of $U_q(sl(2,\C))$. In this
paper, we show that Kashaev's $6j$-symbols are intertwining
operators of local representations of quantum Teichm\"uller
spaces. This relates Kashaev's work with the theory of quantum
Teichm\"uller space, which was developed by Chekhov-Fock, Kashaev
and continued by Bonahon-Liu.

\end{abstract}

\maketitle


\section {Introduction}

Since the early eighties, the theory of quantum invariants of
links and $3$-manifolds has grown up rapidly as a very active
domain of research. However, it is not yet clear which topological
information are carried by these invariants.  One of the most
important conjectural relations between  the topology of a link in
a manifold and its quantum invariants has been constructed by
Kashaev\cite{Kashaev95, Kashaev97} with his Volume Conjecture,
based on a family of complex valued invariants $\{\langle L
\rangle_N\}$ of links $L$ in $S^3$. Later Murakami-Murakami
\cite{MurakamiMurakami01} identified Kashaev's invariants as
evaluations of colored Jones polynomials at the root of unity
$\exp (2\pi i/N)$, which brought the Volume Conjecture to the
forefront. \ignore{
\begin{conj}(Volume Conjecture) Let $L$ be a hyperbolic link in
$S^3$, $J_N(q)$ its $N$-th Jones polynomial, normalized to be $1$
for the unknot. The following holds:
$$\lim_{N\rightarrow \infty} \frac{2\pi }{N}
\log (|J_N(e^{2 \pi i/N})|)=\mathrm{Vol}(S^3\backslash L).$$
\end{conj}
}

The Volume Conjecture is extended to various aspects, such as
Baseilhac-Benedetti \cite{BaseilhacBenedetti04}, Murakami
etc.\cite{Murakami etal02}, Gukov\cite{Gukov05}.  Also see the
references in \cite{Costantino06} for the recent progress on the
Volume Conjecture.

In a different area, the theory of the Quantum Teichm\"uller Space
was developed by Chekhov-Fock \cite{ChekhovFock99} and,
independently, by Kashaev\cite{Kashaev98}. Bonahon-Liu
\cite{BonahonLiu04} investigated the finite dimensional
representation theory of the (exponential version of the) Quantum
Teichm\"uller Space to construct invariants of surface
diffeomorphism; see also \cite{BaiBonahonLiu06}. In computations,
these invariants look very similar to those of
\cite{BaseilhacBenedetti04,Kashaev97,Kashaev98}. Although the
objects involved appear very different, the role played by the
Pentagon Relation and by the quantum dilogarithm function also
hints at a connection between these two points of view. This paper
is devoted to making this connection explicit.

\ignore{ Analog to Turaev-Viro state sum invariants
\cite{TuraevViro92} associated with the quantized universal
enveloping algebra $U_q(sl(2,\C))$,} The Kashaev invariants are
based on the representation theory  of the Weyl algebra $\W$, a
Hopf algebra corresponding to the Borel subalgebra of
$U_q(sl(2,\C))$. In particular, the Hopf algebra structure is used
in a critical way.

 If $\mu: \W\rightarrow \End(V_\mu)$ and $\nu: \W\rightarrow \End(V_\nu)$
 form a \textit{regular} pair of
 (see $\S$\ref{subsec: Weyl algebra} for specific definition)
representations of the Weyl algebra $\W$, then all irreducible
components of the representation $\mu \otimes \nu$ are isomorphic
to a representation  $\mu \nu$. This irreducible component $\mu
\nu$ has  multiplicity $N$, in the sense that the space
$V_{(\mu,\nu)}=\Hom_{\W}(V_{\mu \nu},V_\mu \otimes V_\nu)$ has
dimension $N$. The elements of  $V_{(\mu,\nu)}$ are the
\textit{Clebsch-Gordan operators} of $\mu$ and $\nu$. There is a
natural  \textit{canonical map} $\Omega$
$$\Omega=\Omega(\rho,\mu):V(\rho,\mu)\otimes V_{\rho \mu}
\rightarrow V_\rho \otimes V_\mu
$$defined by $
\Omega(f\otimes v )=f(v)$.

For a regular triple $(\mu,\nu,\sigma)$ of representations of the
Weyl algebra $\mathcal{W}$,  all irreducible components of the
representation $ \mu \otimes \nu \otimes \sigma$ are isomorphic to
a representation $ \mu \nu \sigma$. with multiplicity $N^2$. The
two different ways of grouping terms in $(\mu \otimes \nu) \otimes
\sigma=\mu\otimes (\nu \otimes \sigma)$ lead to two different ways
of grouping together the embedding of $V_{\mu \nu \sigma}$ in
$V_\mu \otimes V_\nu \otimes V_\sigma$. The key ingredient of the
Kashaev invariant of \cite{Kashaev95,BaseilhacBenedetti01} is {\it
Kashaev's $6j$-symbol}  $$R(\mu,\nu,\sigma): V_{(\mu,\nu)}\otimes
V_{(\mu \nu,\sigma)} \rightarrow V_{(\mu,\nu \sigma)}\otimes
V_{(\nu,\sigma)}$$ which describes this correspondence.

More precisely, $R(\mu,\nu,\sigma)$ is defined by the five-term
relation
$$
\Omega_{23}(\nu,\sigma)\Omega_{13}(\mu,\nu
\sigma)R_{12}(\mu,\nu,\sigma)=\Omega_{12}(\mu,\nu)\Omega_{23}(\mu
\nu,\sigma).
$$
involving linear maps $ V_{(\mu,\nu)}\otimes V_{(\mu \nu,\sigma)}
\otimes V_{\mu \nu \sigma} \rightarrow V_\mu \otimes V_\nu \otimes
V_\sigma$, see $\S$\ref{subsec:CGO and 6j-symbols} for details.

 \ignore{
 The
$6j$-symbols naturally satisfy the Pentagon Equation. Kashaev
associated to each  $6j$-symbol a tetrahedron, such that on every
edge there is an irreducible representation of the Weyl algebra.
In this manner, Kashaev constructed the state sum for a given
triangulation of $3$-manifolds $M$, then shown that the state sum
is independent of the triangulation thus gives rising to a well
defined invariant of $M$. }

The quantum Teichm\"uller space $\T_S^q$ of a punctured surface
$S$ is based on completely different data. It is  a
non-commutative  deformation of the algebra of rational functions
on the classical Teichm\"uller space. The non-commutativity is
measured by a parameter $q=e^{\pi i \hbar}$, and the classical
case corresponds to $q=1$.

More precisely, consider an ideal triangulation $\lambda$ of the
surface $S$, with edges $\lambda_1,\cdots,\lambda_n$, and
 triangles $T_1,\cdots, T_m$. The
\textit{Chekhov-Fock algebra}
$\T_{\lambda}^q=\mathbb{C}[X_1,X_2,\ldots,X_n]_{\lambda}^q$ of
$\lambda$ is generated
 by
 $X_1^{\pm1}$, $X_2^{\pm1},\ldots,X_n^{\pm1}$
associated to the edges of $\lambda$, with relations $X_i
X_j=q^{2\sigma_{ij}} X_j X_i$, with integers $\sigma_{ij}\in
\{0,\pm1,\pm2\}$ determined by the combinatorics of the ideal
triangulation. As one moves from one ideal triangulation $\lambda$
to another $\lambda'$, Chekhov and Fock \cite{Fock97,
ChekhovFock99} make this construction independent of the choice of
ideal triangulations by introducing explicit \textit{coordinate
change isomorphisms} $ \Phi_{\lambda \lambda'}^q:\T_{\lambda'}^q
\rightarrow \T_{\lambda}^q$.

 A \textit{local representation} of the
Chekhov-Fock algebra $\T_\lambda^q$ is a certain type of
representation $\rho: \T_\lambda^q
 \rightarrow
\mathrm{End}(V_\lambda) $, where $V_\lambda=V_1\otimes \cdots
\otimes V_m$ and the $V_j$ are associated to the triangles $T_j$;
see $\S$\ref{subsec: intertwining operators} for precise
definition.
 A \textit{local representation} $\rho$ of the quantum Teichm\"uller
space $\T_S^q$ consists of the data of a local representation
$\rho_\lambda: \T_{\lambda}^q \rightarrow \End(V_\lambda)$ for
every ideal triangulation $\lambda$, in such a way that
$\rho_{\lambda}\circ \Phi_{\lambda \lambda'}^q: \T_{\lambda'}^q
\rightarrow \mathrm{End}(V_{\lambda}) $ and $ \rho_{\lambda'}:
\T_{\lambda'}^q \rightarrow \mathrm{End}(V_{\lambda'}) $ are
isomorphic for any pair of ideal triangulations
$\lambda,\lambda'$.
 The \textit{intertwining operator $L_{\lambda
\lambda'}^{\rho}$} of $\rho$ are the  linear maps such that
$$
 L_{\lambda
\lambda'}^{\rho}: V_1\otimes \cdots \otimes V_m \rightarrow
V_1'\otimes \cdots \otimes V_m'$$ that realize the isomorphism
$\rho_{\lambda}\circ \Phi_{\lambda \lambda'}^q \cong
\rho_{\lambda'}$, namely such that $L_{\lambda \lambda'}^{\rho}
\circ \rho_{\lambda}\circ \Phi_{\lambda
\lambda'}^q(X')=\rho_{\lambda'}(X') \circ L_{\lambda
\lambda'}^{\rho} $ for every $X'\in \T_{\lambda'}^q$.

The intertwining operators $L_{\lambda \lambda'}^{\rho} $ are only
defined up to scalar multiplication. They are uniquely determined
by the case when the surface $S$ is a square and when $\lambda$
and $\lambda'$ differ only by the choice of diagonal; see
\cite{BaiBonahonLiu06}.

Our main result is the following theorem:
\begin{thm}
\label{thm: main thm} Let $L:V_1\otimes V_2\rightarrow V_1'
\otimes V_2'$ be a linear map, each vector space with an action of
the Weyl algebra $\W$. Then the following are equivalent:
\begin{enumerate}
\item[1.] there exists a regular pair $(\mu, \nu)$  of irreducible
representations of the Weyl algebra $\W$ and isomorphisms of
$\W$-spaces
$$V_1\cong V_{(\mu,\nu)}, V_2\cong V_{\mu \nu},
V_1'\cong V_\mu, V_2' \cong V_\nu,
$$ for which $L$ corresponds to a scalar multiple of the
canonical map   $$\Omega(\mu,\nu): V_{(\mu,\nu)}\otimes V_{\mu
\nu}  \rightarrow V_\mu \otimes V_\nu.$$
\item[2.] there exists a regular triple $(\mu, \nu, \sigma)$ of
irreducible representations  of the Weyl algebra $\W$ and
isomorphisms of $\W$-spaces
$$V_1\cong V_{(\mu,\nu)}, V_2\cong V_{(\mu \nu, \sigma)},
V_1'\cong V_{(\mu, \nu\sigma)}, V_2' \cong V_{(\nu,\sigma)},
$$
for which $L$ corresponds to a scalar multiple of Kashaev's
$6j$-symbol
$$R(\mu,\nu,\sigma): V(\mu,\nu)\otimes V(\mu \nu,\sigma)
\rightarrow V(\mu,\nu \sigma)\otimes V(\nu,\sigma).$$
\item[3.]  $L$ is an intertwining operator $$L_{\lambda
\lambda'}^\rho: V_1\otimes V_2\rightarrow V_1' \otimes V_2',$$
namely is an isomorphism between $\rho_\lambda \circ \Phi_{\lambda
\lambda'}^q$ and $\rho_{\lambda'}$. Here $\rho=\{ \rho_{\lambda}:
\T_{\lambda}^q \rightarrow \mathrm{End}(V_1\otimes V_2),\,
\rho_{\lambda'}: \T_{\lambda'}^q \rightarrow
\mathrm{End}(V_1'\otimes V_2') \}$ is a local representation
 of the quantum
Teichm\"uller space of the square $S$.
\end{enumerate}
\end{thm}

A  byproduct of Theorem \ref{thm: main thm} is the fact, already
observed by Kashaev\cite{Kashaev95} and Baseilhac
\cite{BaseilhacThesis}, that
 the Clebsch-Gordan operators and the $6j$-symbols
have the same matrix expression in  suitable basis, although their
approach used the representation of the canonical element of the
Heisenberg double of the Weyl algebra.

\medskip
\noindent\textit{Acknowledgments:} The author would like to thank
Francis Bonahon, Xiaobo Liu and Ren Guo  for helpful discussions
during the course of the writing of this paper. The author is very
grateful to the organizers of "A second time around the Volume
Conjecture" for the invitation to talk on this subject, and to the
Department of Mathematics of Louisiana State University for the
kind hospitality.

\section{Kashaev's $6j$-symbol}
\label{sec: Kashaev}

\subsection{The Weyl algebra and its representations}
\label{subsec: Weyl algebra}

The Weyl algebra $\W=\W^q$ is the Hopf algebra over $\C[q,q^{-1}]$
generated by $ X^{\pm 1},Y$ with relation $XY=q^2 YX$. Its
comultiplication, counit and antipode
 are given by
\begin{align*}
&\Delta(X)=X \otimes X, \hspace{4 mm} \Delta(Y)=X^{-1}\otimes Y
+Y\otimes 1; \\
& \epsilon(X)=1, \hspace{4 mm} \epsilon(Y)=0; \\& S(X)=X^{-1},
\hspace{4 mm} S(Y)=-X Y.
\end{align*}
The reader will notice that the notation for the Weyl algebra $\W$
 is slightly different from that in
\cite{BaseilhacBenedetti01, Kashaev95}, but it is clearly
equivalent and better suited to our purpose.

 A \textit{representation} of the
Hopf algebra $\mathcal{W}$ is a left $\mathcal{W}$-module, here
$\mathcal{W}$ is simply regarded  as an algebra. \ignore{The terms
subrepresentation, irreducible representation, complete
reducibility, etc., have their usual meanings.} Given a
representation $\mu: \W\rightarrow \End(V_\mu)$ of $\W$, the
action of $W\in \W$ on an element $v\in V_\mu$ is denoted by
$W_\mu \cdot v=\mu(W)(v)$, or simply $W\cdot v$ if no ambiguity on
the representation.

 The comultiplication  of $\mathcal{W}$ allows one to
make the tensor product of representations of $\mathcal{W}$ into a
new representation. If $\mu: \W\rightarrow \End(V_\mu)$ and $\nu:
\W\rightarrow \End(V_\nu)$ are representations of $\mathcal{W}$,
then the representation $\mu\otimes \nu$ is defined by
$$
W\cdot (v\otimes w)=\Delta(W)\cdot (v\otimes w)= \sum W_{(1)}\cdot
v \otimes W_{(2)}\cdot w
$$
for any $W\in \W$, $v\in V_\mu$ and $w\in V_\nu$. Here we use
Sweedler's sigma notation
$$ \Delta(W) = \sum W_{(1)}\otimes W_{(2)}. $$
For general properties of representation theory of Hopf algebras,
we refer readers to Chari-Pressley\cite{ChariPressley94} or
Montgomery\cite{Montgomery93}.

 Any linear space $V$ can be given a \textit{trivial}
representation of  $\mathcal{W}$  by
$$ W\cdot v=\epsilon(W)v  \hsp \forall W\in \W, v\in V.$$
 It is of interest to consider non-trivial finite dimensional
representations of $\W$. They exist only when $q^2$ is a root of
 unity. We consequently assume that  $q^2$ be a primitive $N$-th root of
 unity for some positive integer $N$.

\ignore{ Not every representation of the Weyl algebra
$\mathcal{W}$ is complete reducible. Such an example could be
$$V=
\mathrm{Span}_\mathbb{C}\{ v, Y\cdot v, \cdots, Y^{N-1}\cdot v
\}$$
 with conditions $X\cdot v=q^2 v$, $Y^{N-1}\cdot v\ne 0$ and $Y^N\cdot v=0$. It is
easy to see that $\mathcal{W}$ acts trivially on  the subspace
spanned by a single vector $ Y^{N-1}\cdot v$, but there exists no
$\mathcal{W}$-invariant compliment subspace in $V$ for this
one-dimensional subspace. In the example above, the action of $Y$
on the space is not invertible. This observation leads to the
following definition \cite{Kashaev95,BaseilhacThesis}. }

A representation $\mu: \W \rightarrow \mathrm{End}(V_{\mu})$ over
a finite dimensional complex vector space $V_{\mu}$ is  \textit
{cyclic} if the operators $\mu(X)$ and $\mu(Y)$ are invertible.
 A sequence
of cyclic representations $(\mu_1,\cdots, \mu_n)$ is \textit
{regular} if for any  $1\leq i <n$ and $1\leq j \leq n-i$ the
representations $\mu_i\otimes \mu_{i+1}\otimes \cdots\otimes
\mu_{i+j}$ is cyclic.

\ignore{ A representation $\mu: \W\rightarrow \End(\C^N)$ is
\textit{standard} with parameters $(x,y)$ if there exist $x,y \in
\mathbb{C}^*$ such that
$$
X_\mu=x A,\hspace{6 mm} Y_\mu=y B
$$
 where
}

The following is completely elementary.

\begin{prop} \label{prop:Rep Weyl algebra}
 Let $q^2$ be a primitive $N$-th root of unity.
\begin{enumerate}
 \item[1.]
 Every cyclic representation of the Weyl
algebra $\W$ is completely reducible, and its irreducible factors
are cyclic.
 \item[2.]
 Up to isomorphism, an irreducible cyclic representation $\mu: \W\rightarrow \End(V_\mu)$ of the
Weyl algebra
 is completely determined by the complex numbers $x_\mu,y_\mu \in \mathbb{C}^*$  such that
$$
X^N_\mu=x_\mu \cdot \mathrm{Id}_{V_\mu}, \hspace{3 mm}
Y^N_\mu=y_\mu \cdot \mathrm{Id}_{V_\mu}.
$$
 \item[3.] More precisely, the irreducible cyclic representation
 $\mu$ is isomorphic to a representation $\mu': \W\rightarrow
 \End(\C^N)$ defined by
 $$X_{\mu'}=x_\mu^{1/N}A, \hsp Y_{\mu'}=y_\mu^{1/N}B$$
 where $x_\mu^{1/N}$ and $y_\mu^{1/N}$ are arbitrary N-th roots of
 $x_\mu$ and $y_\mu$, and where $A, B \in \End(\C^N)$ are the  unitary matrices
defined by $$A e_k=q^{2k}e_k,\,\, B e_k=e_{k+1} \hspace{2
mm}\mathrm{for} \hspace{2 mm} 0\le k \le N-1 $$ on the standard
basis $\{e_0, \ldots, e_{N-1}\}$ of $\C^N$, with indices
considered modulo $N$.
\end{enumerate}
\qed
\end{prop}

\subsection{The Clebsch-Gordan operators and $6j$-symbols}
\label{subsec:CGO and 6j-symbols}

Suppose that two irreducible representations $\mu$ and $\nu$ form
a regular pair, and are classified by parameters $(x_\mu,y_\mu)$
and $(x_\nu,y_\nu)$ as in Proposition \ref{prop:Rep Weyl algebra}.
  From the calculation of
$$
\Delta(X^N)=X^N \otimes X^N, \hspace{4 mm}
\Delta(Y^N)=X^{-N}\otimes Y^N +Y^N\otimes 1,
$$
we see that all the irreducible components of $\mu\otimes \nu$ are
isomorphic to the representation $\mu \nu$ classified by the
parameters
$$
x_{\mu \nu}=x_\mu x_\nu, \hspace{3 mm} y_{\mu \nu}=x_\mu ^{-1}
y_\nu +y_\mu.
$$
 Thus
\begin{lem}
 If $(\mu,\nu)$ is a regular pair of irreducible
representations of the Weyl algebra $\W$, then the representation
$\mu \otimes \nu: \W \rightarrow \End(V_{\mu} \otimes V_{\nu})$
splits as a direct sum of $N$ representations isomorphic to the
representation $\mu \nu$.

As a consequence, the space $V_{(\mu,\nu)}=\Hom_\W(V_{\mu \nu}, \,
V_\mu \otimes V_\nu)$ of $\W$-equivariant linear maps  $V_{\mu
\nu} \rightarrow V_\mu \otimes V_\nu$ is a vector space of
dimension $N$. \qed
\end{lem}

The elements of $V_{(\mu,\nu)}$ are called the
\textit{Clebsch-Gordan operators}. \ignore{the elements in
$\Hom_\W(V_\mu \otimes V_\nu, \,  V_{\mu \nu})$ (which are
projectors) are called \textit{dual Clebsch-Gordan operators}.}
 The linear map
$\Omega=\Omega(\mu,\nu):V_{(\mu,\nu)}\otimes V_{\mu \nu}
\rightarrow V_\mu \otimes V_\nu $ given by
$$
\Omega(f\otimes v )=f(v), \hsp \forall f\in V_{(\mu,\nu)}, v\in
V_{\mu \nu}
$$ is  the \textit{canonical map}.

 The Kashaev $6j$-symbol is associated to  a regular
sequence $(\mu,\nu,\sigma)$ of three irreducible representations
of the Weyl algebra $\mathcal{W}$.  In this situation, the diagram
\begin{align*}
& \nearrow V_{\mu \nu}\otimes V_\sigma
   \rightarrow V_\mu \otimes V_\nu \otimes
V_\sigma\\
V_{\mu \nu  \sigma} & \hspace{1.5 in} \parallel \\
&\searrow
 V_\mu \otimes V_{\nu \sigma}
\rightarrow V_\mu \otimes V_\nu \otimes V_\sigma
\end{align*}
provides natural isomorphisms between $\Hom_{\W}(V_{\mu \nu
\sigma}, V_\mu \otimes V_\nu \otimes V_\sigma)$ and both
$V_{(\mu,\nu)}\otimes V_{(\mu \nu,\sigma)}$ and $V_{(\mu,\nu
\sigma)}\otimes V_{(\nu,\sigma)}$.
 The {\it
Kashaev's $6j$-symbol} of $(\mu,\nu,\sigma)$ is defined as the
resulting isomorphism
$$R(\mu,\nu,\sigma): V_{(\mu,\nu)}\otimes V_{(\mu \nu,\sigma)} \rightarrow
V_{(\mu,\nu \sigma)}\otimes V_{(\nu,\sigma)}.$$ In other words,
$R(\mu,\nu,\sigma)$ is a linear operator defined  by  equation
$$
\Omega_{23}(\nu,\sigma)\Omega_{13}(\mu,\nu
\sigma)R_{12}(\mu,\nu,\sigma)=\Omega_{12}(\mu,\nu)\Omega_{23}(\mu
\nu,\sigma),
$$
which is better
 illustrated as the commutativity of
pentagon diagram
\begin{align*}
& \stackrel{\Omega_{23}(\nu, \sigma)}{\swarrow} V_\mu \otimes
V_{(\nu, \sigma)}\otimes V_{\nu \sigma}
   \stackrel{\Omega_{13}(\mu,\nu \sigma)}{\longleftarrow} V_{(\mu,\nu \sigma)}\otimes V_{(\nu,\sigma)}\otimes
V_{\mu \nu
\sigma}\\
V_\mu \otimes V_\nu \otimes V_\sigma & \hspace{2.7 in} \uparrow R_{12}(\mu,\nu,\sigma) \\
&\stackrel{\Omega_{12}(\mu,\nu)}{\nwarrow}
 V_{(\mu,\nu)}\otimes V_{\mu \nu}\otimes V_\sigma
\stackrel{\Omega_{23}(\mu \nu, \sigma)}{\longleftarrow}
V_{(\mu,\nu)}\otimes V_{(\mu \nu, \sigma)}\otimes V_{\mu \nu
\sigma}
\end{align*}
Here we use the standard notation that,  if an operator $F$ acts
on the $i$-th and $j$-th components ($1\le i <j \le 3$) of the
tensor product of spaces $U\otimes V \otimes W$, then $F_{ij}$
denotes  the operator of $U\otimes V \otimes W$ that acts by $F$
on the $i$-th and $j$-th components, and by the identity on the
remaining component of the tensor product.

\ignore{If every space of Clebsch-Gordan operators in the diagram
is given a trivial representation of the Weyl algebra, then all
the arrows become $\W$-isomorphisms, according to Lemma \ref{lem:
Omega}. }

\section{The quantum Teichm\"uller space}
\label{sec: Chekhov-Fock}

\subsection{The triangle algebra and its representations}
\label{subsec:triangle algebra}

The \textit{triangle algebra $\mathcal{T}$} is the algebra over
$\mathbb{C}$ generated by $X^{\pm 1},Y^{\pm 1},Z^{\pm 1}$ with
relations
$$XY=q^2YX,\,\, YZ=q^2ZY,\,\, ZX=q^2XZ.
$$
 The center of the triangle algebra $\T$ is generated by the
\textit{principal cental element } $H=q^{-1}XYZ$ and the three
other elements $X^N,Y^N$ and $Z^N$.

\ignore{If we compare two algebras $\mathcal{W}$ and
$\mathcal{T}$, the generator $D$ in $\mathcal{W}$ is not
invertible, while all the generators of triangle algebra are
invertible. We will obtain similar results on the representation
theory of the triangle algebra $\mathcal{T}$.}

The finite dimensional representation theory of the quantum
Teichm\"uller space works best when $q^N=(-1)^{N+1}$, which is
slightly more restrictive than the condition that $q^2$ is a
primitive $N$-th root of unity
\cite{BaiBonahonLiu06,BonahonLiu04}. Consequently we henceforth
assume that $q$ is a primitive $N$-th root of $(-1)^{N+1}$.

Let $\{e_0, \ldots, e_{N-1}\}$ be the standard basis of  $\C^N$,
with indices modulo $N$. Let $A, B, C \in \End(\C^N)$ be the
unitary matrices given by
$$A e_k=q^{2k}e_k, B e_k=e_{k+1}, C e_k=q^{1-2k}e_{k-1}.$$
Note that $AB=q^2BA, BC=q^2CB, CA=q^2AC$  and
$$A^N=B^N=C^N=q^{-1}ABC=\Id_{\C^N}.$$

\ignore{ Assume that $\rho: \mathcal{T} \rightarrow
\End(V_{\rho})$ is an irreducible representation of the triangle
algebra $\mathcal{T}$. Let us pick up $v \in V_{\rho}$, any one of
the eigenvectors of $X$ with eigenvalue $x$,  then
$$\rho(X_i)=x_iA,\hspace{2 mm} \rho(X_{i+1})=x_{i+1}B, \hspace{2
mm} \rho(X_{i+2})=x_{i+2}C $$ under the basis
$$\left\{v_i,\frac{X_{i+1}}{x_{i+1}}(v_i),\left(\frac{X_{i+1}}{x_{i+1}}\right)^2(v_i),
\cdots,\left(\frac{X_{i+1}}{x_{i+1}}\right)^{N-1}(v_i) \right\}.$$
}

The following analogue of Proposition \ref{prop:Rep Weyl algebra}
is elementary.

\begin{prop}
\label{prop:Rep triangle algebra} Let $q^2$ be a primitive $N$-th
root of unity.

\begin{enumerate}
\item[1.]
 Up to isomorphism, an  irreducible representation
 $\rho: \T\rightarrow \End(V_\rho)$ of the
 triangle algebra $\T$ is completely determined by
the complex numbers $x_\rho,y_\rho, z_\rho, h_\rho\in \C^*$ with
relation $h_\rho^N=x_\rho y_\rho z_\rho$, such that
$$\hspace{15 mm}
X^N_\rho=x_\rho\cdot \mathrm{Id}_{V_\rho},\hsp Y^N_\rho=y_\rho
\cdot \mathrm{Id}_{V_\rho},\hsp Z^N_\rho=z_\rho \cdot
\mathrm{Id}_{V_\rho}, \hsp H_\rho=h_\rho \cdot
\mathrm{Id}_{V_\rho}.
$$
In addition, any number $x,y, z, h\in \C^*$ with $h^N=xyz$ can be
realized by such a representation.
 \item[2.] More precisely, the
irreducible  representation
 $\rho$ is isomorphic to a representation $\rho': \T\rightarrow
 \End(\C^N)$ defined by
$$ \hspace{15 mm} X_{\rho'}=x_\rho^{1/N}A, \hsp Y_{\rho'}=y_\rho^{1/N}B,
\hsp Z_{\rho'}=z_\rho^{1/N}C,
 \hsp H_{\rho'}=h_\rho \cdot \mathrm{Id}_{\C^N}
 $$
 where $x_\rho^{1/N}, y_\rho^{1/N}$ and $z_\rho^{1/N}$ are arbitrary N-th roots of
 $x_\rho,y_\rho$ and $z_\rho$ satisfying $x_\rho^{1/N}\cdot y_\rho^{1/N}\cdot z_\rho^{1/N}=h_\rho$,
  and where $A, B, C \in \End(\C^N)$ are the
unitary matrices defined as above.
\end{enumerate} \qed
\end{prop}
 The number $h_\rho$ is called the \textit{central
charge} of the representation $\rho$.

\ignore{ We associate the triangle algebra with a triangle, whose
edges are indexed by $\lambda_1$, $\lambda_2$ and $\lambda_3$
clockwise. Then an irreducible representation $\rho:\T \rightarrow
\End(V_\rho)$ is completely determined by the complex number
weights $z_i$ on the edge $\lambda_i$ ($i=1,2,3$) and by the
central charge $h$ which is an $N$-th root of $z_1 z_2 z_3$. }

\subsection{The quantum Teichm\"uller space}
\label{subsec: QTS}

 Let $S$ be an oriented surface of genus
$g$ with $p\ge 1$ punctures, obtained by removing $p$ points
$\{v_1, v_2, \dots, v_p\}$ from the closed oriented surface
$\overline{S}$ of genus $g$. \ignore{The \textit{Teichm\"uller
space $\mathcal{T}(S)$} is the space of isotopy classes of
complete hyperbolic metrics on $S$.}An \textit{ideal
triangulation} $\lambda$ of $S$ consists of finitely many disjoint
simple arcs $\lambda_1,\lambda_2,\dots,\lambda_n$ going from
puncture to puncture and decomposing $S$ into triangles. Any ideal
triangulation has precisely $n=-3\chi(S)=6g+3p-6$ arcs, where
$\chi(S)$ is the Euler characteristic of the punctured surface.
Let $\Lambda(S)$  denote the set of isotopy classes of all ideal
triangulations.

Following the
terminology of \cite{BonahonLiu04, Liu04}, the
\textit{Chekhov-Fock algebra} associated to the ideal
triangulation $\lambda$ is the algebra
$\T_\lambda^q=\mathbb{C}[X_1,X_2,\ldots,X_n]_{\lambda}^q$ over
$\mathbb{C}$
 generated
 by
 $X_1^{\pm1},X_2^{\pm1},\ldots,X_n^{\pm1}$ respectively
associated to the edges of $\lambda$ with relations $X_i
X_j=q^{2\sigma_{ij}} X_j X_i$, where $\sigma_{ij}\in
\{0,\pm1,\pm2\}$ are integers determined by the combinatorics of
the ideal triangulation. This algebra has a well-defined fraction
division algebra $\mathbb{C}(X_1,X_2,\dots,X_n)_{\lambda}^q$,
consisting of all formal rational fractions in variables $X_i$
that skew-commute according to the relations $X_i X_j =
q^{2\sigma_{ij}} X_j X_i$.

 As one moves from one ideal triangulation
$\lambda$ to another $\lambda'$, Chekhov and Fock \cite{Fock97,
ChekhovFock99} (as developed in \cite{Liu04}) introduce explicit
\textit{coordinate change isomorphisms} $ \Phi_{\lambda
\lambda'}^q:\mathbb{C}(X_1',X_2',\dots,X_n')_{\lambda'}^q
\rightarrow \mathbb{C}(X_1,X_2,\dots,X_n)_{\lambda}^q. $ These are
algebra isomorphisms which satisfy the natural property that
$\Phi_{\lambda  \lambda''}^q= \Phi_{\lambda \lambda'}^q
 \circ \Phi_{\lambda' \lambda''}^q $ for any ideal triangulations
$\lambda, \lambda', \lambda''\in \Lambda(S)$.
 These algebra isomorphisms are
essentially unique \cite{Bai05}, once we require them to satisfy a
certain number of natural conditions.

The \textit{quantum Teichm\"{u}ller space} of $S$ is the algebra
defined in a triangulation independent way as
\begin{align*}
\mathcal{T}_{S}^q &= \left(\bigsqcup_{\lambda
          \in \Lambda(S) }
         \mathbb{C}(X_1,X_2,\dots,X_n)_{\lambda}^q
   \right) / \sim \\
&= \{ (\lambda,X): \lambda \in \Lambda(S), \,\, X \in
\mathbb{C}(X_1,X_2,\dots,X_n)_{\lambda}^q
   \}/ \sim
\end{align*}
where the equivalence relation $\sim$ is defined by
$$
(\lambda, X) \sim (\lambda', X') \Leftrightarrow X=\Phi_{\lambda
\lambda'}^q(X').
$$

\subsection{Local representations  of the quantum
Teichm\"uller space}
\label{subsec: intertwining operators}
 Every ideal triangulation $\lambda$ of $S$ has exactly $m=-2\chi(S)$
triangles. Each triangle $T_j$ determines a triangle algebra
 $\mathcal{T}_j$, with generators associated to the three sides of
$T_j$. The Chekhov-Fock algebra has a natural embedding
$\iota_\lambda$ into the tensor product of triangle algebras
$$\iota_{\lambda}: \T_{\lambda}^q
\rightarrow \bigotimes_{j=1}^m \mathcal{T}_j,$$ see \cite{Bai06,
BaiBonahonLiu06}.  To describe $\iota_{\lambda}$, consider the
generator $X_i\in \T_{\lambda}^q $ corresponding to the edge
$\lambda_i$ of $\lambda$.  There are two possible cases:
\begin{enumerate}
\item
 The edge $\lambda_i$ separates two triangles $T_j$ and $T_{j'}$.
 Define
$ \iota_{\lambda}(X_i)= X_{ij} \otimes X_{ij'} $ where $X_{ij}(\,
X_{ij'} \,\, \mathrm {resp.})$ is the generator of $\T_j(\,
\T_{j'} \,\, \mathrm {resp.})$ associated to the edge $\lambda_i$.
\item The two sides of $\lambda_i$ belong to the same
  triangle $T_j$.
Let  $X_{ij_1}$ and $X_{ij_2}$ be the generators of
$\mathcal{T}_{j}$ associated to
  those two sides, indexed such that $X_{ij_1} X_{ij_2}=q^2 X_{ij_2} X_{ij_1}$,
then define $ \iota_{\lambda}(X_i)= q^{-1}X_{ij_1}X_{ij_2} =q
X_{ij_2} X_{ij_1}$
\end{enumerate}
By convention, when describing an element $Z_1 \otimes  \dots
\otimes Z_m$ of $\T_{1} \otimes  \dots \otimes \T_{m}$, we omit in
the tensor product those $Z_j$ that are equal to the identity
element $1$ of $\T_{j}$.

Suppose that we are given an irreducible representation $\rho_j:
\mathcal{T}_j\rightarrow \mathrm{End}(V_j)$ for every triangle
$T_j$. The tensor product $\otimes_{j=1}^m \rho_j$ restricts to a
representation $ \rho: \T_{\lambda}^q \rightarrow \End(V_\lambda)
$ where $V_\lambda=V_1\otimes \cdots \otimes V_m$. \ignore{ which
satisfies
$$\rho(q^{-\sum_{i<j} \sigma_{ij}} X_1 X_2\cdots X_n)=h
\mathrm{Id}.$$ }
 By definition, a \textit{local representation} of
$\T_{\lambda}^q$ is any representation obtained in this way.

 A \textit{local representation} of the quantum Teichm\"uller
space $\T_S^q$ consists of the data of a local representation
$\rho_\lambda: \T_{\lambda}^q \rightarrow \End(V_\lambda)$ of the
Chekhov-Fock algebra for every ideal triangulation $\lambda$, in
such a way that $\rho_{\lambda}\circ \Phi_{\lambda \lambda'}^q:
\T_{\lambda'}^q \rightarrow \mathrm{End}(V_{\lambda}) $ and $
\rho_{\lambda'}: \T_{\lambda'}^q \rightarrow
\mathrm{End}(V_{\lambda'}) $ are isomorphic for any pair of ideal
triangulations $\lambda,\lambda'$.

\ignore{
 Note that special care should
be taken to make sense of $\rho_{\lambda}\circ \Phi_{\lambda
\lambda'}^q: \T_{\lambda'}^q \rightarrow \mathrm{End}(V_{\lambda})
$ over finite dimensional space $V_{\lambda}$. Readers are
referred to \cite{BaiBonahonLiu06} for rigorous treatments. }

 Let $\rho=\left\{ \rho_\lambda:  \T_{\lambda}^q \rightarrow
\End(V_\lambda)\right\}_{\lambda\in \Lambda(S)}$ be a local
representation of $\T_S^q$. For every  pair of ideal
triangulations $\lambda,\lambda'$, \ignore{ the local
representations of distinct Chekhov-Fock algebras $
\rho_{\lambda}: \T_{\lambda}^q \rightarrow \mathrm{End}(V_1\otimes
\cdots \otimes V_m) $ and $ \rho_{\lambda'}: \T_{\lambda'}^q
\rightarrow \mathrm{End}(V_1'\otimes \cdots \otimes V_m') $
  are related in a way such that $\rho_{\lambda}\circ
\Phi_{\lambda \lambda'}^q $ and $\rho_{\lambda'}$ of
$\T_{\lambda'}^q$ are isomorphic.
}
 an
\textit{intertwining operator $L_{\lambda \lambda'}^{\rho}$} for
the local representation $\rho$ is defined as the linear map
$$
 L_{\lambda
\lambda'}^{\rho}: V_1\otimes \cdots \otimes V_m \rightarrow
V_1'\otimes \cdots \otimes V_m'$$ such that $L_{\lambda
\lambda'}^{\rho} \circ \rho_{\lambda}\circ \Phi_{\lambda
\lambda'}^q(X')=\rho_{\lambda'}(X') \circ L_{\lambda
\lambda'}^{\rho} $ for every $X'\in \T_{\lambda'}^q$. This
$L_{\lambda \lambda'}^{\rho}$ is actually what is denoted by
$L_{\lambda \lambda'}^{\rho \rho}$ in \cite{BaiBonahonLiu06}.

\ignore{The intertwining operators can be determined only up to
scalar multiplication. We will consequently write $A\doteq  B$ to
say that two linear maps $A$ and $B$ are equal up to
multiplication by a scalar.}

Theorem $20$ of \cite{BaiBonahonLiu06} shows that  the family of
intertwining operators $L_{\lambda\lambda'}^\rho$ is uniquely
determined, provided that we impose certain natural conditions on
them. In particular, the intertwining operators can be explicitly
obtained, for all surfaces $S$, from the intertwining operator
associated to the case where $S$ is a square and $\lambda$,
$\lambda'$ are ideal triangulations of $S$ differing from each
other by different choices of a diagonal of that square.

Here, a \textit{square} is the surface obtained from a closed disk
by removing $4$ punctures from its boundary.

 \ignore{ Up to scalar multiplication, there exists a unique
family of intertwining operators $L_{\lambda\lambda'}^\rho$,
indexed by local representations $\rho$ of quantum Teichm\"uller
spaces $\T_S^q$ of surfaces $S$ and by ideal triangulations
$\lambda$, $\lambda'$ of the same surface $S$, such that:
\begin{enumerate}
\item \emph{(Composition Relation)} for any three ideal
triangulations $\lambda$, $\lambda'$, $\lambda''$ of the same
surface $S$ and for any local representation $\rho$ of $\T_S^q$,
$$ L_{\lambda\lambda''}^\rho \dot =L_{\lambda'\lambda''}^\rho \circ L_{\lambda\lambda'}^\rho;$$
\item \emph{(Fusion Relation)}  if $S$ is obtained from another
surface $S_0$ by fusing $S_0$ along certain components of
$\partial S_0$, if the representation $\rho$ of the quantum
Teichm\"uller space $\T_S^q$ is obtained by fusing a
representation $\rho_0$ of $\T_{S_0}^q$, and if $\lambda$,
$\lambda'$ are two ideal triangulations of $S$ obtained by fusing
ideal triangulations $\lambda_0$, $\lambda_0'$ of $S_0$, then
$L_{\lambda\lambda'}^\rho \dot= L_{\lambda_0\lambda_0'}^{\rho_0}$.
\end{enumerate}
}

\subsection{The quantum Teichm\"uller space of the square}
\label{subsec: QTS of square} We consequently analyze the case of
the square in more details.
\begin{figure}[h]
\begin{center}
\includegraphics[width=5 in]{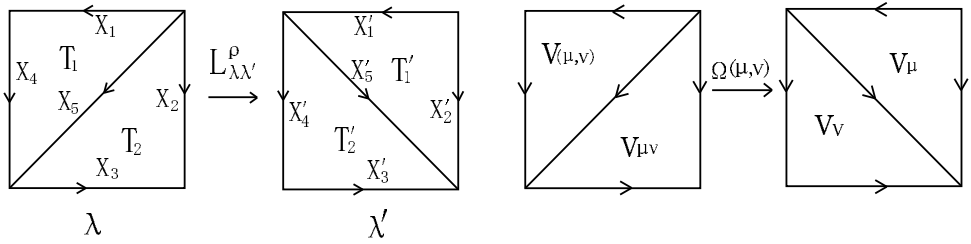}
\caption{The intertwining operator and the canonical map}
\label{fig:IO and Omega}
\end{center}
\end{figure}

Let $\lambda$ and $\lambda'$ be the two ideal triangulations of
the square $S$ represented in Figure \ref{fig:IO and Omega}, with
the indexing of edges and triangles indicated there. Let
$\rho_\lambda: \T_{\lambda}^q\rightarrow \End(V_1\otimes V_2)$ and
$\rho_{\lambda'}: \T_{\lambda'}^q\rightarrow \End(V_1'\otimes
V_2')$ be local representations of the Chekhov-Fock algebra of
$\lambda$ and $\lambda'$.

The \textit{principal cental element} of the Chekhov-Fock algebra
$\T_{\lambda}^q$ is $H=q^{-\sum_{i<j}\sigma_{ij}}X_1 X_2\cdots
X_5=X_1 X_2\cdots X_5$. Similarly, the principal cental element of
Chekhov-Fock algebra $\T_{\lambda'}^q$ is
$H'=q^{-\sum_{i<j}\sigma_{ij}}X_1' X_2'\cdots X_5'= q^{-2}X_1'
X_2'\cdots X_5'$. Compare with the principal cental element of the
triangle algebra in $\S$\ref{subsec:triangle algebra}.

The following is proved in \cite{BaiBonahonLiu06}, although the
proof is elementary.

\begin{prop}
\label{prop: classify local rep} Up to isomorphism, the local
representation $\rho_\lambda:\T_{\lambda}^q\rightarrow
\End(V_1\otimes V_2)$ is uniquely determined by the numbers $x_1,
x_2, x_3,x_4, h_{\rho_\lambda} \in \C^*$ such that
$\rho_\lambda(X_i^N)=x_i\cdot \Id \,\, (i=1,2,3,4)$ and
$\rho_\lambda(H)=h_{\rho_\lambda}\cdot \Id$. In addition, any five
tuple of numbers $x_1, x_2, x_3,x_4, h \in \C^*$ are associated to
a local representation in this way. \qed
\end{prop}
The  number $h_{\rho_\lambda}$ is called the \textit{central
charge} of the local representation $\rho_\lambda$. Note that
$\rho_\lambda(X_5^N)=x_5\cdot \Id$ where
$x_5=h_{\rho_\lambda}^N/(x_1 x_2 x_3 x_4)$.

The same result holds for $\lambda'$, so that $\rho_{\lambda'}$ is
classified by numbers $x_1', x_2', x_3',x_4', h_{\rho_{\lambda'}}
\in \C^*$.

\begin{prop}
\label{prop: condition of parameters}
 $\rho_\lambda \circ
\Phi_{\lambda \lambda'}^q$ is isomorphic to $\rho_{\lambda'}$ if
and only if
$$\left\{
\begin{array}{ll}
x_1' & =  (1+x_5)x_1 \\
x_2' & =  (1+x_5^{-1})^{-1} x_2 \\
x_3' & =  (1+x_5)x_3 \\
x_4' & =  (1+x_5^{-1})^{-1} x_4 \\
 h_{\rho_{\lambda'}} & =
h_{\rho_\lambda}
\end{array}
\right.$$ where $x_5=h_{\rho_\lambda}^N/(x_1 x_2 x_3 x_4)$.
\end{prop}

\begin{proof} This is an immediate consequence of Theorem $12$ of
\cite{BaiBonahonLiu06}, and of the specific form of the coordinate
change isomorphism $\Phi_{\lambda \lambda'}^q$ for the square.
\end{proof}

We should point out that the role of $x_i$ and $x_i'$ are
symmetric, in the sense that the above formula is equivalent to
the following one
$$\left\{
\begin{array}{ll}
x_1 & = (1+{x'}_5^{-1})^{-1}x_1' \\
x_2 & =  (1+x_5') x_2' \\
x_3 & =  (1+{x'}_5^{-1})^{-1} x_3' \\
x_4 & =  (1+x_5') x_4' \\
h_{\rho_{\lambda}} & =  h_{\rho_{\lambda'}}
\end{array}
\right.$$ where $x_5'=h_{\rho_{\lambda'}}^N/(x_1' x_2' x_3'
x_4')=x_5^{-1}$.

\begin{lem}
\label{lem: h=1}

Suppose that $\rho_\lambda \circ \Phi_{\lambda \lambda'}^q$ is
isomorphic to $\rho_{\lambda'}$ by an intertwining operator
$L_{\lambda \lambda'}^q: V_1\otimes V_2 \rightarrow V_1' \otimes
V_2'$. Then there exist unique local representations
$\rho_\lambda': \T_{\lambda}^q\rightarrow \End(V_1\otimes V_2)$
and $\rho_{\lambda'}': \T_{\lambda'}^q\rightarrow \End(V_1'\otimes
V_2')$ with central charges equal to $1$ and such that
$$\rho_\lambda'(X_i)=\rho_\lambda(X_i), \hsp \rho_{\lambda'}'(X_i')=\rho_{\lambda'}(X_i')$$
for $i=1,2,3,4$.

In addition, the same $L_{\lambda \lambda'}^q$ defines an
isomorphism between $\rho_\lambda' \circ \Phi_{\lambda
\lambda'}^q$ and $\rho_{\lambda'}'$.
\end{lem}
\begin{proof}
Define $\rho_\lambda'$ by the property that $\rho_\lambda'(X_i)$
is equal to $\rho_\lambda(X_i)$ if $1\le i \le 4$, and
$\rho_\lambda'(X_5)= \rho_\lambda(X_5)/h$. Similarly, define
$\rho_{\lambda'}'$.
 \end{proof}

As a consequence, we can restrict our attention to local
representations with central charge $1$.

\section{The Weyl algebra and the triangle algebra}
\label{sec: Weyl to triangle}

As an algebra, the Weyl algebra $\W$ is clearly very similar to
the triangle algebra $\T$. Let us make this correspondence
explicit.

Let $T$ be a triangle, with an additional structure of an arrow on
each edge arranged so that the following holds: the three vertices
of the triangle can be ordered so that the ordering goes
counterclockwise, and each edge is oriented from the lower vertex
to the higher vertex. Figure \ref{fig:IO and Omega} offers many
examples.

With this additional data, any irreducible regular representation
 $\mu: \W\rightarrow \End(V_\mu)$  of the Weyl algebra gives a
unique representation $ \T_T \rightarrow \End(V_\mu)$, which we
will also denote by $\mu$, of the triangle  algebra $\T_T$
associated to the triangle $T$. Indeed,  the generators $X,Y,Z$ of
the triangle algebra $\T_T$ can be arranged such that $X$ is
associated to the edge going from the lowest vertex to the middle
one, while $Y$ is associated to the edge going from lowest vertex
to the highest one.

In this setting, define
 the representation $\mu: \T_T
\rightarrow \End(V_\mu)$ (induced from the representation $\mu: \W
\rightarrow \End(V_\mu)$) by
$$X_\mu=X_\mu, \,\, Y_\mu=Y_\mu, \,\, Z_\mu=q Y_\mu^{-1} X_\mu^{-1}.$$
Note that the central charge of $\mu$ is $1$, as
$H_\mu=(q^{-1}XYZ)_\mu=\mathrm{Id}_{V_\mu}$.  Every irreducible
representation of $\T_T$ with central charge $1$ is obtained in
the way.

\section {An action of the Weyl algebra on the Clebsch-Gordan operators}
\label{sec: Weyl on CGO}

Assume that $\mu:\W\rightarrow \End(V_\mu)$, $\nu:\W\rightarrow
\End(V_\nu)$ form a regular pair of irreducible representations of
the Weyl algebra. Let $V_{(\mu,\nu)}=\Hom_{\W}(V_{\mu \nu}, V_\mu
\otimes V_\nu)$ be the corresponding space of Clebsch-Gordan
operators.

Consider the square $S$, with the two triangulations $\lambda$,
 $\lambda'$ and the edge orientations indicated in Figure \ref{fig:IO and Omega}.
 Associate to each triangle of $\lambda$,
 $\lambda'$ vector spaces $V_\mu$, $ V_\nu$, $V_{\mu \nu}$
 and $V_{(\mu,\nu)}$ as in Figure \ref{fig:IO and Omega}.

 Note that the spaces $V_\mu$, $ V_\nu$, $V_{\mu \nu}$  come with
 an action of the Weyl algebra $\W$, and therefore, using the
 arrows indicated, with an action of the triangle algebra of the
 corresponding triangle.

 If, in addition, the space of Clebsch-Gordan
operators $V_{(\mu,\nu)}$ is also endowed with an action of  the
Weyl algebra $\W$, this will define two local representations
$\rho_\lambda:\T_{\lambda}^q \rightarrow \End(V_{(\mu,\nu)}\otimes
V_{\mu\nu})$ and $\rho_\lambda':\T_{\lambda'}^q\rightarrow
\End(V_\mu \otimes V_\nu)$ of corresponding Chekhov-Fock algebras.

\begin{prop}
\label{prop: Omega as IO} Assume that $\mu:\W\rightarrow
\End(V_\mu)$ and $\nu:\W\rightarrow \End(V_\nu)$ form a regular
pair of irreducible representations of the Weyl algebra. Then
there exists a unique action the Weyl algebra on the space of
Clebsch-Gordan operators $V_{(\mu,\nu)}$, such that the canonical
map $\Omega(\mu,\nu):V_{(\mu, \nu)}\otimes V_{\mu \nu} \rightarrow
V_\mu \otimes V_\nu $ is an intertwining operator for the local
representation $\rho$ of the quantum Teichm\"uller space of the
square, which consists of $\rho_\lambda:\T_{\lambda}^q \rightarrow
\End(V_{(\mu,\nu)}\otimes V_{\mu\nu})$ and
$\rho_\lambda':\T_{\lambda'}^q\rightarrow \End(V_\mu \otimes
V_\nu)$.
\end{prop}

\begin{proof} Let us denote the possible action of the Weyl algebra on the
Clebsch-Gordan operators by  $(\mu,\nu): \W \rightarrow
\End{V_{(\mu,\nu)}}$. If the canonical map
$\Omega=\Omega(\mu,\nu)$ is an intertwining operator for the local
representation $\rho$, then by definition $\Omega \circ
\rho_\lambda \circ \Phi_{\lambda \lambda'}^q= \rho_\lambda' \circ
\Omega$, or equivalently, $\Omega \circ \rho_\lambda =
\rho_\lambda' \circ (\Phi_{\lambda \lambda'}^q)^{-1} \circ
\Omega$.

A typical set of generators of $\T_{\lambda}^q$ consists of $X_1,
\cdots, X_5$ as described in $\S$\ref{subsec: QTS of square}. This
set corresponds to $\{X\otimes 1, 1\otimes Y, 1\otimes Z, Z\otimes
1, Y\otimes X\}$ in the algebra tensor $\T_{T_1}\otimes \T_{T_2}$
via the natural embedding of $\S$\ref{subsec: intertwining
operators}, with the arrows as in Figure \ref{fig:IO and Omega}.

 In particular, for any $f\in V_{(\mu,\nu)}$ and
$v\in V_{\mu \nu}$, the equation
$$\Omega \circ \rho_\lambda (X\otimes 1)(f\otimes v)=
\rho_\lambda' \circ (\Phi_{\lambda \lambda'}^q)^{-1} (X\otimes
1)\circ \Omega(f\otimes v)$$ gives that $$(X_{(\mu,\nu)}\cdot f)
(v)=(X^{-1}_\mu \otimes \mathrm{Id}_{V_\nu} + Y_\mu\otimes
Y^{-1}_\nu)^{-1} \cdot f(v),$$ and the equation $$\Omega \circ
\rho_\lambda (Y\otimes X)(f\otimes v)= \rho_\lambda' \circ
(\Phi_{\lambda \lambda'}^q)^{-1} (Y\otimes X)\circ \Omega(f\otimes
v)$$ gives that $$ (Y_{(\mu,\nu)}\cdot f) (X_{\mu \nu}\cdot
v)=(Z_\mu^{-1} \otimes Y_\nu^{-1} )\cdot f(v).$$ If we require
that $Y_{(\mu,\nu)}\cdot f$ is still in $V_{(\mu,\nu)}$, namely is
$\W$-equivariant, then it forces that
$$ (Y_{(\mu,\nu)}\cdot f) (
v)=(X_\mu^{-1} \otimes X_\nu^{-1})\circ (Z_\mu^{-1} \otimes
Y_\nu^{-1} )\cdot f(v)=(q^{-1}Y_\mu \otimes X_\nu^{-1}
Y_\nu^{-1})\cdot f(v).$$

Hence the possible action $(\mu,\nu)$ is completely determined by
the formulas
\begin{align*}(X_{(\mu,\nu)}\cdot f) (v)&=(X^{-1}_\mu \otimes
\mathrm{Id}_{V_\nu} + Y_\mu\otimes
Y^{-1}_\nu)^{-1} \cdot f(v), \\
 (Y_{(\mu,\nu)}\cdot f) (v)&=(q^{-1}Y_\mu \otimes
X_\nu^{-1} Y_\nu^{-1})\cdot f(v).
\end{align*}

 The condition
that $(\mu, \nu)$ forms a regular pair of irreducible
representations of the Weyl algebra is critical, which guarantees
that both the operators $X^{-1}_\mu \otimes \mathrm{Id}_{V_\nu} +
Y_\mu\otimes Y^{-1}_\nu$ and $q^{-1}Y_\mu \otimes X_\nu^{-1}
Y_\nu^{-1}$ are invertible in $\End(V_\mu \otimes V_\nu)$.

 We need to show that the above formula really defines an action of the Weyl
algebra  on the space of Clebsch-Gordan operators $V_{(\mu,\nu)}$,
namely that for any $f\in V_{(\mu,\nu)}=\Hom_{\W}(V_{\mu \nu},
V_\mu \otimes V_\nu)$, these operators $X_{(\mu,\nu)}\cdot f$ and
$Y_{(\mu,\nu)}\cdot f$ are still in $V_{(\mu,\nu)}$, or,  are
$\W$-equivariant.

The comultiplication $\Delta$ of the Weyl algebra $\W$ has the
property that $\Delta(X)=X \otimes X$ and $\Delta(Y)=X^{-1}\otimes
Y +Y\otimes 1 $, see $\S$ \ref{subsec: Weyl algebra}.  It is
straightforward to check that both the elements $X^{-1} \otimes 1
+ Y\otimes Y^{-1}$ and $q^{-1}Y \otimes X^{-1} Y^{-1}$ commute
with either $\Delta(X)$ or $\Delta(Y)$ in $\W \otimes \W$,
therefore commute with $\Delta(W)$ for any $W\in \W$.

Now $X_{(\mu,\nu)}\cdot f \in V_{(\mu,\nu)}$ is immediate from the
following calculation
\begin{align*} (X_{(\mu,\nu)}\cdot f) (a_{\mu \nu}\cdot
v)&= (X^{-1}_\mu \otimes \mathrm{Id} + Y_\mu\otimes
Y^{-1}_\nu)^{-1} \cdot f(a_{\mu \nu}\cdot  v)\\
&=(X^{-1}_\mu \otimes \mathrm{Id} + Y_\mu\otimes Y^{-1}_\nu)^{-1}
\circ
\Delta(a)\cdot f(v)\\
&=\Delta(a) \circ (X^{-1}_\mu \otimes \mathrm{Id} + Y_\mu\otimes
Y^{-1}_\nu)^{-1} \cdot f(v)\\
&=\Delta(a) \circ (X_{(\mu,\nu)}\cdot f) (v).
\end{align*}
Same argument proves that $Y_{(\mu,\nu)}\cdot f \in
V_{(\mu,\nu)}$.

By establishing the representation of the Weyl algebra on the
space of Clebsch-Gordan operators, we have shown that $$\Omega
\circ \rho_\lambda = \rho_\lambda' \circ (\Phi_{\lambda
\lambda'}^q)^{-1} \circ \Omega$$  holds for $X\otimes 1$ and
$Y\otimes X$. The reader can easily check  that the equation holds
for other generators of the Chekhov-Fock algebra $\T_{\lambda}^q$.

This completes the proof the proposition.
\end{proof}

Calculation shows that $$ X_{(\mu,\nu)}^N= (x_\mu^{-1}+y_\mu
y_\nu^{-1})^{-1} \cdot \Id, \hsp Y_{(\mu,\nu)}^N= y_\mu x_\nu^{-1}
y_\nu^{-1}\cdot \Id.$$ That means, the representation $(\mu,\nu)$
of the Weyl algebra is classified by the parameter $ (
x_{(\mu,\nu)}= (x_\mu^{-1}+y_\mu  y_\nu^{-1})^{-1}, y_{(\mu,\nu)}=
y_\mu x_\nu^{-1} y_\nu^{-1})$ in the sense of Proposition
\ref{prop:Rep Weyl algebra}. Once again, we note that the  number
$x_\mu^{-1}+y_\mu y_\nu^{-1}\ne 0$  due to the fact that $\mu$ and
$\nu$ forms a regular pair.

\section{The canonical maps, Kashaev's $6j$-symbols, and intertwining operators}
\label{sec: Omega 6j and  IO}

\begin{thm}
Let $L:V_1\otimes V_2\rightarrow V_1' \otimes V_2'$ be a linear
map, each vector space with an action of the Weyl algebra $\W$.
Then the following are equivalent:
\begin{enumerate}
\item[1.] there exists a regular pair $(\mu, \nu)$  of irreducible
representations of the Weyl algebra $\W$ and isomorphisms of
$\W$-spaces
$$V_1\cong V_{(\mu,\nu)}, V_2\cong V_{\mu \nu},
V_1'\cong V_\mu, V_2' \cong V_\nu,
$$ for which $L$ corresponds to a scalar multiple of the
canonical map   $$\Omega(\mu,\nu): V_{(\mu,\nu)}\otimes V_{\mu
\nu}  \rightarrow V_\mu \otimes V_\nu.$$
\item[2.] there exists a regular triple $(\mu, \nu, \sigma)$ of
irreducible representations  of the Weyl algebra $\W$ and
isomorphisms of $\W$-spaces
$$V_1\cong V_{(\mu,\nu)}, V_2\cong V_{(\mu \nu, \sigma)},
V_1'\cong V_{(\mu, \nu\sigma)}, V_2' \cong V_{(\nu,\sigma)},
$$
for which $L$ corresponds to a scalar multiple of Kashaev's
$6j$-symbol
$$R(\mu,\nu,\sigma): V(\mu,\nu)\otimes V(\mu \nu,\sigma)
\rightarrow V(\mu,\nu \sigma)\otimes V(\nu,\sigma).$$
\item[3.]  $L$ is an intertwining operator $$L_{\lambda
\lambda'}^\rho: V_1\otimes V_2\rightarrow V_1' \otimes V_2',$$
namely is an isomorphism between $\rho_\lambda \circ \Phi_{\lambda
\lambda'}^q$ and $\rho_{\lambda'}$. Here $\rho=\{ \rho_{\lambda}:
\T_{\lambda}^q \rightarrow \mathrm{End}(V_1\otimes V_2),\,
\rho_{\lambda'}: \T_{\lambda'}^q \rightarrow
\mathrm{End}(V_1'\otimes V_2') \}$ is a local representation
 of the quantum
Teichm\"uller space of the square $S$.
\end{enumerate}
\end{thm}

\begin{proof} [Proof of $1\Rightarrow 3$] This is a consequence of
Proposition \ref{prop: Omega as IO}.
\end{proof}

\begin{proof} [Proof of $3\Rightarrow 1$] Let us pick
$V_\nu=V_1'$, $V_\nu=V_2'$.

First, we need to show that there exist isomorphisms of
$\W$-spaces between $V_{\mu \nu}$ and $V_2$, also between
$V_{(\mu,\nu)}$ and $V_1$. Let us consider the relationship
between the corresponding parameters.

Notice that in every representation of the triangle algebra which
is induced from an irreducible regular representation of the Weyl
algebra, the product of the three parameters $x,y,z$ is equal to
$1$, see $\S$ \ref{sec: Weyl to triangle}. Thus from $V_\nu=V_1'$,
$V_\nu=V_2'$,  the following holds
$$x_\mu= x_1', \,\, y_\mu= x_2', \,\, x_\nu=x_4', \,\, y_\nu={x'}_3^{-1} {x'}_4^{-1}.$$
Now combining with the relation
$$
x_1  = (1+{x'}_5^{-1})^{-1}x_1', \, x_2 =  (1+x_5') x_2', \, x_3 =
(1+{x'}_5^{-1})^{-1} x_3', \, x_4  =  (1+x_5') x_4'
$$ from $\S$ \ref{subsec: QTS of square} (where $x_5'=1/(x_1' x_2' x_3' x_4')$),
the relation $$ x_{(\mu,\nu)}= (x_\mu^{-1}+y_\mu y_\nu^{-1})^{-1},
y_{(\mu,\nu)}= y_\mu x_\nu^{-1} y_\nu^{-1}$$ from $\S$\ref{sec:
Weyl on CGO} and the relation
$$
x_{\mu \nu}=x_\mu x_\nu, \hspace{3 mm} y_{\mu \nu}=x_\mu ^{-1}
y_\nu +y_\mu
$$ from $\S$\ref{subsec:CGO and 6j-symbols}, we immediately get
$$x_{(\mu,\nu)}=x_1, \,\, y_{(\mu,\nu)}=x_1^{-1}x_4^{-1}, \,\,
x_{\mu \nu}=x_2^{-1} x_3^{-1}, \,\, y_{\mu \nu}=x_2.$$ By
Proposition \ref{prop:Rep Weyl algebra}, this is exactly the
relationship among the parameters which gives the desirable
isomorphisms of $\W$-spaces.

Then, we need to show that $L$ corresponds to the canonical map
$\Omega(\mu,\nu)$, up to scalar multiplication.

Proposition \ref{prop: Omega as IO} tells us that
$\Omega(\mu,\nu)$ is an intertwining operator for the local
representation $\rho$ of the quantum Teichm\"uller space of the
square, which consists of $\rho_\lambda:\T_{\lambda}^q \rightarrow
\End(V_{(\mu,\nu)}\otimes V_{\mu\nu})$ and
$\rho_\lambda':\T_{\lambda'}^q\rightarrow \End(V_\mu \otimes
V_\nu)$. But any irreducible representation of the quantum
Teichm\"uller space of the square is of dimension $N^2$
\cite{BonahonLiu04}, hence both $\rho_\lambda \circ \Phi_{\lambda
\lambda'}^q$ and $\rho_\lambda'$ are irreducible. Then Schur's
lemma implies that the intertwining operator is unique, up to
scalar multiplication.  This finishes our proof.
\end{proof}

\begin{proof}[Proof of $2 \Rightarrow 3$]
 Assume that  $(\mu, \nu, \sigma)$ forms a
regular triple
 of irreducible representations  of the Weyl
algebra $\W$. We interpret the defining equation of  Kashaev's
$6j$-symbol
$$R(\mu,\nu,\sigma): V(\mu,\nu)\otimes V(\mu \nu,\sigma)
\rightarrow V(\mu,\nu \sigma)\otimes V(\nu,\sigma)$$ as the
commutativity of the Figure \ref{fig:6j}.
\begin{figure}[h]
\begin{center}
\includegraphics[width=4 in]{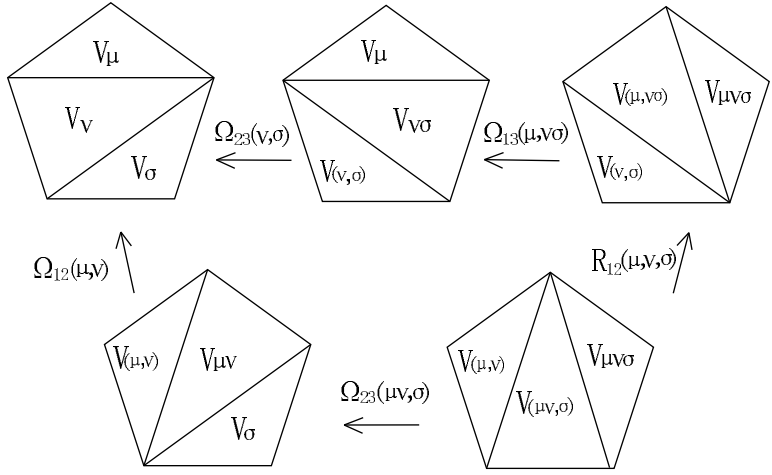}
\caption{$6j$-symbol as intertwining operator} \label{fig:6j}
\end{center}
\end{figure}

By the result of $1 \Rightarrow 3$, every $\Omega$-map is an
intertwining operator between the local representations of the
Chekhov-Fock algebras associated to the five distinct ideal
triangulations of the pentagon in Figure \ref{fig:6j}.

Note that, by the property of the coordinate change isomorphisms
$\Phi_{\lambda' \lambda}^q = \left(\Phi_{\lambda \lambda'}^q
\right)^{-1}$,  the inverse of an intertwining operator is still
an intertwining operator, since $L_{\lambda \lambda'}^\rho \circ
\rho_{\lambda} \circ \Phi_{\lambda
\lambda'}^q=\rho_{\lambda'}\circ L_{\lambda \lambda'}^\rho$
implies $\left(L_{\lambda \lambda'}^\rho \right)^{-1} \circ
\rho_{\lambda'} \circ \Phi_{\lambda' \lambda}^q = \rho_{\lambda}
\circ \left(L_{\lambda \lambda'}^\rho \right)^{-1} $.

Then  $$R(\mu,\nu,\sigma)\otimes \Id_{V_{\mu \nu
\sigma}}=\Omega_{13}^{-1}(\mu,\nu \sigma) \circ
\Omega_{23}^{-1}(\nu,\sigma)\circ \Omega_{12}(\mu,\nu) \circ
\Omega_{23}(\mu \nu,\sigma)$$ is the composition of the
intertwining operators, so again is an  intertwining operator
between the local representations of corresponding Chekhov-Fock
algebras of the pentagon. Consequently, Kashaev's $6j$-symbol
$R(\mu,\nu,\sigma)$ is an  intertwining operator for a local
representation of the quantum Teichm\"uller space of the square.
\end{proof}

\begin{proof}[Proof of $3 \Rightarrow 2$]
Let $\rho=\{ \rho_{\lambda}, \rho_{\lambda'} \}$ be a local
representation of the quantum Teichm\"uller space of the square
$S$. Assume that the local representation $\rho_{\lambda'}:
\T_{\lambda'}^q \rightarrow \mathrm{End}(V_1'\otimes V_2')$ is
classified by the number $x_1', x_2',x_3',x_4' \in \C^*$
($h_{\rho_{\lambda'}}=1$) in the sense of Proposition \ref{prop:
classify local rep}.  If there exists a  regular triple $(\mu,
\nu, \sigma)$ of irreducible representations  of the Weyl algebra
$\W$, such that the representation $(\mu,\nu \sigma)\otimes
(\nu,\sigma)$ of the Weyl algebra induces the local representation
$\rho_{\lambda'}$ of the Chekhov-Fock algebra, then by Proposition
\ref{prop: condition of parameters} the parameters necessarily
have relation
$$\left\{
\begin{array}{lll}
x_1' & = x_{(\mu,\nu\sigma)} &= \left( x_\mu^{-1}+y_{\mu}(x_\nu^{-1}y_\sigma +y_\nu)^{-1}\right)^{-1} \\
x_2' & = y_{(\mu,\nu\sigma)} &= y_{\mu}\cdot x_\nu^{-1} x_\sigma^{-1} \cdot (x_\nu^{-1}y_{\sigma}+y_\nu)^{-1} \\
x_3' & = x_{(\nu,\sigma)} &= \left(x_\nu^{-1}+y_{\nu}y_\sigma^{-1} \right)^{-1} \\
x_4' & = z_{(\nu,\sigma)} &=(x_\nu^{-1}+y_\nu y_\sigma^{-1})\cdot
y_\nu^{-1} x_\sigma y_\sigma
\end{array}
\right. $$ The second equality of each row is derived from the
relations from $\S$\ref{sec: Weyl on CGO} and
 $\S$\ref{subsec:CGO and 6j-symbols}, see the proof of $3 \Rightarrow 1$.

This is easily solved as
 $$\left\{
\begin{array}{ll}
x_\mu  = {x'}_1^{-1}-x_2' x_3' x_4' ({x'}_3^{-1} y_\nu^{-1}
y_\sigma -1)^{-1},
&y_\mu= x_2'  x_4' ( {x'}_3^{-1}y_\nu^{-1} - y_\sigma^{-1})^{-1}  \\
x_\nu  = ({x'}_3^{-1}-y_\nu y_\sigma^{-1})^{-1}, &y_\nu= y_\nu\\
x_\sigma  = x_3' x_4' y_\nu y_\sigma^{-1},  &y_\sigma= y_\sigma \\
\end{array}
\right. $$ for every $y_\nu, y_\sigma$ such that $$y_\nu \cdot
y_\sigma \cdot (y_\sigma y_\nu^{-1} - x_3')\cdot (y_\sigma
y_\nu^{-1}-(x_1' x_2' x_3' x_4'+1)x_3')\ne 0.$$

Any such solution  to these equations defines regular irreducible
representations $\mu, \nu,\sigma$ of the Weyl algebra and
$\W$-isomorphisms $$V_1\cong V_{(\mu,\nu)}, V_2\cong V_{(\mu \nu,
\sigma)}, V_1'\cong V_{(\mu, \nu\sigma)}, V_2' \cong
V_{(\nu,\sigma)}.
$$
By the proof of $2\Rightarrow 3$, we know that Kashaev's
$6j$-symbol $R(\mu,\nu,\sigma): V(\mu,\nu)\otimes V(\mu
\nu,\sigma) \rightarrow V(\mu,\nu \sigma)\otimes V(\nu,\sigma)$ is
an intertwining operator, and we already see that, in the case of
the square, the intertwining operator is unique up to scalar
multiplication.  It follows that this $6j$-symbol
$R(\mu,\nu,\sigma)$ corresponds to a scalar multiple of $L$.
\end{proof}

\end{document}